\documentclass[12pt]{article}

\def\chi{{\mathcal{X}}}
\def\calh{{\mathcal{H}}}
\def\calv{{\mathcal{V}}}
\def\1{^{-1}}

\def\wedge{{\Lambda}}
\def\calb{{\mathcal{B}}}
\def\cald{{\mathcal{D}}}

\def\cala{\mathcal{A}}
\def\({\left(}
\def\){\right)}
\def\pf{\n{\bf Proof.} }
\def\vsp{\vspace*{1,5mm}\\ }
\def\vspp{\vspace*{2mm}\\ }

\def\bk{\bigskip }
\def\mk{\medskip }
\def\sk{\smallskip }

\def\n{\noindent }
\def\dd{\displaystyle}

\def\D{{\Delta}}

\def\barr{\begin{array}}
\def\earr{\end{array}}

\def\bit{\begin{itemize}}
\def\itemi{\item[{\rm(i)}]}
\def\itemii{\item[{\rm(ii)}]}
\def\itemiii{\item[{\rm(iii)}]}
\def\itemiv{\item[{\rm(iv)}]}

\def\eit{\end{itemize}}

\def\D{{\Delta}}

\usepackage{amsmath, amsthm, amscd, amsfonts, amssymb}

\newtheorem{theorem}{Theorem}[section]

\newtheorem{corollary}[theorem]{Corollary}

\theoremstyle{definition}
\newtheorem{definition}[theorem]{Definition}
\newtheorem{example}[theorem]{Example}

\newtheorem{remark}[theorem]{Remark}

\def\1{^{-1}}
\def\vsp{\vspace*{2mm}\\ }

\def\calf{{\mathcal{F}}}
\def\calo{{\mathcal{O}}}

\def\E{{\mathbb{E}}}
\def\rr{{\mathbb{R}}}
\def\nn{{\mathbb{N}}}
\def\9{{\infty}}
\def\G{{\Gamma}}
\def\lbb{{\lambda}}
\def\a{{\alpha}}

\def\g{{\gamma}}
\def\wt{\widetilde}
\def\ov{\overline}
\def\vf{{\varphi}}
\def\oo{{\omega}}
\def\ooo{{\Omega}}
\def\pp{{\partial}}
\def\D{{\Delta}}

\def\barr{\begin{array}}
\def\earr{\end{array}}

\def\dd{\displaystyle}
\def\bk{\bigskip }
\def\sk{\smallskip}
\def\n{\noindent }

\def\pas{\mathbb{P}\mbox{-a.s.}}

\def\vsp{\vspace*{2mm}\\ }
\def\ff{\forall }
\def\({\left(}
\def\){\right)}
\def\<{\left<}
\def\>{\right>}

\title{A  splitting algorithm for~stochastic partial differential equations driven by linear multiplicative noise}
\author{Viorel Barbu\thanks{Octav Mayer Institute of
Mathematics (Romanian Academy) and Al.I. Cuza University of Ia\c si, Romania; e-mail: vbarbu41@gmail.com} \and Michael R\"ockner\thanks{Fakult\"at f\"ur Mathematik, Universit\"at Bielefeld, D-33501 Bielefeld, Germany; e-mail: roeckner@math.uni-bielefeld.de}}
\date{}

\begin{document}
\maketitle
\begin{abstract}
We study the convergence of a Douglas-Rachford type splitting algorithm for the infinite dimensional stochastic differential equation $$dX+A(t)(X)dt=X\,dW\mbox{ in }(0,T);\  X(0)=x,$$ where $A(t):V\to V'$ is a nonlinear, monotone, coercive and demicontinuous operator with sublinear growth and $V$ is a real Hilbert space with the dual $V'$. $V$ is densely and continuously embedded in the Hilbert space $H$ and  $W$ is an $H$-valued Wiener process. The general case of a maximal monotone operators $A(t):H\to H$ is also investigated.\sk
\end{abstract}

\section{Introduction}\label{s1}
We consider here the stochastic differential equation
\begin{equation}\label{e1.1}
\barr{l}
dX(t)+A(t)X(t)dt=X(t)dW(t),\ t\in(0,T),\vsp
X(0)=x,\earr\end{equation}in a real separable Hilbert space $H$, whose elements are functions or distributions on a bounded and open set $\calo\subset\rr^d$ with smooth boundary $\pp\calo$. In particular, $H$ can be any of the spaces $L^2(\calo)$, $H^1_0(\calo),$ $H\1(\calo)$, $H^1(\calo)$, $k=1,2,...,$ with the corresponding Hilbertian structure. Here $H^1_0(\calo),$ $H^k(\calo)$ are the standard  $L^2$-Sobolev spaces on $\calo$, and $W$ is a Wiener process of the form
\begin{equation}\label{e1.2}
W(t,\xi)=\sum^\9_{j=1}\mu_je_j(\xi)\beta_j(t),\ \xi\in\calo,\ t\ge0,\end{equation}
where $\{\beta_j\}^\9_{j=1}$ is an independent system of real-valued Brownian motions on a probability space $\{\ooo,\calf,\mathbb{P}\}$ with natural filtration $(\calf_t)_{t\ge0}$. Here, $e_j\in C^2(\ov\calo)\cap H$, $j\in\mathbb{N}$, is an orthonormal basis in $H$, and $\mu_j\in\rr,\ j=1,2,...$.

The following hypotheses will be in effect throughout this work.

\bit\itemi There is  a Hilbert space $V$ with dual $V'$ such that $V\subset H$, continuously and densely. Hence $V\subset H\ (\equiv H')\subset V'$  continuously and densely.
\itemii $A:[0,T]\times V\times\ooo\to V'$ is {\it progressively measurable}, i.e., for every $t\in[0,T]$, this operator restricted to $[0,t]\times V\times\ooo$ is $\calb([0,t])\otimes\calb(V)\otimes\calf_t$-measurable.
    \itemiii There is $\delta\ge0$ such that, for each $t\in[0,T]$, $\oo\in\ooo$, the operator $u\mapsto\delta u+A(t,\oo)u$ is monotone and {\it demicontinuous} (that is, strongly-weakly continuous) from $V$ to $V'$.

Moreover, there are $\a_i,\g_i\in\rr,\ i=1,2,3,$ $\a_1>0$, such that, $\pas$,
\begin{eqnarray}
\<A(t,\oo)u,u\>&\ge&\a_1|u|^2_V+\a_2|u|^2_H+\a_3,\ \ff u\in V,\ t\in[0,T],\quad\label{e1.3}\\[2mm]
|A(t,\oo)u|_{V'}&\le&\g_1|u|_V+\g_2,\hspace*{19mm} \ff u\in V,\ t\in[0,T].\quad\label{e1.4}
\end{eqnarray}
\itemiv $e^{\pm W(t)}$ is, for each $t$, a multiplier in $V$ and a   multiplier in $H$ such that   there exists an $(\calf_t)$-adapted, $\rr_+$-valued process  $Z(t)$, $t\in[0,T]$, with $\mathbb{E}\left[\sup\limits_{t\in[0,T]}|Z(t)| \right]<\9$ for all $r\in[1,\9)$ and such that, $\pas$,
\begin{equation}
\barr{lcll}
|e^{\pm W(t)}y|_V&\le &Z(t)|y|_V,& \ff t\in[0,T],\ \ff y\in V, \vsp
|e^{\pm W(t)}y|_H&\le &Z(t)|y|_H,&\ff t\in[0,T],\ \ff y\in H.\earr\label{e1.5}
\end{equation}
\eit
One assumes also that, for each $\oo\in\ooo$, the function $t\to e^{\pm W(t)}$ is $H$-valued continuous on $[0,T]$.

Throughout in the following, $|\cdot|_V$ and $|\cdot|_{V'}$ denote the norms of $V$ and $V'$, respectively, and by $\<\cdot,\cdot\>$ we denote the duality pairing between $V$ and $V'$ with $H$ as pivot space; on $H\times H$, $\<\cdot,\cdot\>$ is just the scalar product of $H$. The norms of $V$ and $V'$  are denoted by $|\cdot|_H$ and $|\cdot|_V$, $|\cdot|_{V'}$, respectively, $\calb(H)$, $\calb(V)$ etc. are the classes of Borel sets in the corresponding spaces.

As regards the orthonormal basis $\{e_j\}^\9_{j=1}$ in \eqref{e1.2}, we assume that there exist $\wt\g_j\in[1,\9)$ such that
\begin{equation}\label{e1.7}
|ye_j|_H\le\wt\g_j |y|_H,\ \ff y\in H,\ j=1,2,...,\ \nu:=\sum^\9_{j=1}\mu^2_j\wt\g^2_j|e_j|^2_\9<\9.
\end{equation}
and  we assume also that
\begin{equation}\label{e1.8}
\mu:=\frac12\sum^\9_{j=1}\mu^2_je^2_j \end{equation}
 is a multiplier in $V$, $V'$  and  $H$.

It should be  noted that $X\,dW=\sigma(X)d\wt W$ where $\sigma:H\to L_2(H)$ (the space of Hilbert-Schmidt operators on $H$) is defined by
$$\sigma(u)v=\sum^\9_{j=1}\mu_j u\<v,e_j\>e_j,\ \ff v\in H,$$and so, $\wt W=\sum\limits^{\9}_{j=1}e_j\beta_j$ is a cylindrical Wiener process on $H$ (see \cite{4a}).

\begin{definition}\label{d1.1} \rm By a {\it solution} to \eqref{e1.1} for $x\in H$, we mean an $(\calf_t)_{t\ge0}$-adapted process $X:[0,T]\to H$ with continuous sample paths which satisfies
\begin{eqnarray}
&X(t)+\int^t_0A(s)X(s)ds=x+\int^t_0X(s)dW(s),\ t\in[0,T],\label{e1.10}\\[2mm]
&  X\in L^2((0,T)\times\ooo;V).\label{e1.11}
\end{eqnarray}
\end{definition}

The stochastic integral arising in \eqref{e1.10} is considered in It\^o's sense.

In \cite{2}, the authors developed an operatorial approach to \eqref{e1.1} under the more general hypotheses than (i)--(iv) above.
As a special case (see Theorem  3.1 in \cite{2}), we have

\begin{theorem}\label{t1.2} Under Hypotheses {\rm(i)--(iv)}, for each $x\in H$, equation \eqref{e1.1} has a unique solution $X$ $($in the sense of Definition {\rm\ref{d1.1}$)$.} Moreover, the function $t\mapsto e^{-W(t)}X(t)$ is $V'$-absolutely continuous on $[0,T]$ and
\begin{equation}\label{e1.12}
\E\int^T_0\left|e^{W(t)}\,\frac d{dt}\,\(e^{-W(t)}X(t)\)\right|^2_{V'}dt<\9.\end{equation}
\end{theorem}

In a few words, the method developed in \cite{2} is the following. By the transformation
\begin{equation}\label{e1.13}
X(t)=e^{W(t)}y(t),\ t\ge0,\end{equation}
one reduces equation \eqref{e1.1} to the random differential equation
\begin{equation}\label{e1.14}
\barr{l}
\dd\frac{dy}{dt}\,(t)+e^{-W(t)}A(t)\(e^{W(t)}y(t)\)+\mu y(t)=0,\mbox{\ a.e. }t\in(0,T),\vsp
y(0)=x,\earr\end{equation}
and treat \eqref{e1.14} as an operatorial equation of the form
\begin{equation}\label{e1.15}
\calb y+\cala y=0\end{equation}
in a suitable Hilbert space $\calh$ of stochastic processes on $[0,T]$. Here, $\cala$ and $\calb $ are maximal monotone operators suitable defined from $\calv$ to $\calv'$, where $(\calv,\calv')$ is a dual pair of spaces such that $\calv\subset\calh\subset\calv'$ with dense and continuous embeddings.

The operatorial form \eqref{e1.15} of equation \eqref{e1.14} suggests to approximate the solution $y$  by the Douglas--Rachford splitting algorithm (\cite{4}--\cite{6}).

The exact form and   convergence of the corresponding splitting algorithm for equation \eqref{e1.15} will be given below in Section 2. As seen later on in Theorem \ref{t2.1}, it leads to a convergent splitting algorithm for the stochastic differential equation \eqref{e1.1}.

In this way, the operator theoretic approach to equation \eqref{e1.1} written in  the form \eqref{e1.15} allows to design a convergent splitting scheme for equation \eqref{e1.1} inspired by the Rockafellar \cite{7} proximal point algorithm for nonlinear operatorial equations (on these lines see also \cite{3}). By our knowledge, the splitting algorithm  obtained here for the stochastic equation is new and might have implications in numerical approximation of stochastic PDEs.

\bk\n{\bf Notations.} If $U$ is a Banach space, we denote by $L^p(0,T;U),$  $1\le p\le\9$,\break the space of all $L^p$-integrable $U$-valued functions on $(0,T)$. The space\break $L^p((0,T)\times\ooo;U)$ is defined similarly.   We refer to \cite{1} for notation and standard results of the theory of maximal monotone operators in Banach spaces. If $\calo$ is an open domain of $\rr^d$, we denote by $W^{1,p}(\calo)$, $1\le p\le\9$ and $H^1(\calo),H^{-1}(\calo)$ the standard Sobolev spaces on $\calo$.

\section{Main results}
\setcounter{equation}{0}
\setcounter{theorem}{0}

Without loss of generality, we may assume that, besides assumptions (i)--(iii), $A(t)$ satisfies also the strong monotonicity condition
\begin{equation}\label{e2.1}
\<A(t)u-A(t)v,u-v\>\ge\nu|u-v|^2_H,\ \ff u,v\in V,\end{equation}
 where $\nu>0$ is given by \eqref{e1.7}. (In fact, as easily seen, by the substitution   $X\to\exp(-(\nu+\delta)t)X$ with a suitable $\delta$, equation \eqref{e1.1} can be rewritten as$$dX+\wt A(t)X\,dt=X\,dW,$$where the operator $X\to \wt A(t)X=e^{-(\nu+\delta)t}A(t)(e^{(\nu+\delta)t}X)
 +(\nu+\delta)X$ satisfies conditions  (i)--(iii) and \eqref{e2.1}.)

 We associate with equation \eqref{e1.1} the following splitting algorithm
\begin{equation}\label{e2.2}
\barr{l}
\lbb dZ_{n+1}+J(Z_{n+1})dt+\lbb \nu Z_{n+1}dt
=\lbb Z_{n+1}dW-\lbb A(t)X_ndt\vsp
\hfill+\lbb\nu X_n\,dt+J(X_n)dt,\ t\in(0,T),\vsp
Z_{n+1}(0)=x,\ n=0,1,...\earr\end{equation}
\begin{equation}\label{e2.3}
\barr{l}
\lbb A(t)X_{n+1}(t)+J(X_{n+1}(t))-\lbb\nu X_{n+1}(t)\vsp\qquad
=J(Z_{n+1}(t))+\lbb A(t)X_n(t)-\lbb\nu X_n(t),\earr\end{equation}where  $X_0\in L^2((0,T)\times\ooo;V)$ is   $(\calf_t)_{t\ge0}$-adapted and arbitrary. Here, the parameter $\lbb>0$ is arbitrary but fixed and $J:V\to V'$ is the canonical isomorphism of the space $V$ onto its dual $V'$.

Taking into account assumptions (i)--(iii) and \eqref{e2.1}, which, in particular, implies  that the operator $\G_0:L^2(0,T;V)\to L^2(0,T;V')$, $\G_0 u=\lbb A(t)u+J(u)-\lbb\nu u,$  $u\in L^2(0,T;V)$, is demicontinuous, locally bounded, and with inverse continuous, we see that the sequence $(Z_n,X_n)$ is well defined by \eqref{e2.2}, \eqref{e2.3} and we   have also
\begin{equation}\label{e2.4}
X_n,Z_n\in L^2((0,T)\times\ooo;V)\mbox{ and }Z_n\in L^2(\ooo;C([0,T];H)),\ n=1,2,...\end{equation}
Moreover, the processes $X_n,Z_n$ are $(\calf_t)_{t\ge0}$-adapted on $[0,T]$.

Theorem \ref{t2.1} is the main result.

\begin{theorem}\label{t2.1} Under Hypotheses {\rm(i)--(iv)} and \eqref{e2.1}, assume that $x\in V$ and $\lbb>0$. If $(X_n,Z_n)$ is the sequence defined by \eqref{e2.2}, \eqref{e2.3}, we have for $n\to\9$
\begin{equation}\label{e2.5}
X_n\to X\mbox{ weakly in }L^2((0,T)\times\ooo;V),\end{equation}where $X$ is the solution to equation \eqref{e1.1} given by Theorem {\rm\ref{t1.2}.} Assume further that the operator $u\to A(t)u$ is odd, that is, $A(t)(-u)=-A(t)u$, $\ff u\in V.$ Then, for $n\to\9$,
\begin{equation}\label{e2.6}
X_n\to X\mbox{ strongly in }L^2((0,T)\times\ooo;V).\end{equation}
\end{theorem}

The splitting scheme \eqref{e2.2}--\eqref{e2.3} reduces the approximation of problem \eqref{e1.1} to a sequence of simpler linear equations. In fact, at each step $n$, one should solve a linear stochastic differential equation of the form
\begin{equation}\label{e2.7}
\barr{l}
dZ_{n+1}+\dd\frac1\lbb\,J(Z_{n+1})dt+\nu Z_{n+1}dt=Z_{n+1}dW+F_ndt,\ t\in(0,T),\vsp
Z_{n+1}(0)=x,\earr\end{equation}and the stationary random equation \eqref{e2.3}, where $$F_n=-\lbb A(t)X_n+\lbb\nu X_n+J(X_n).$$

By It\^o's formula (see, e.g.,  \cite{2}), equation \eqref{e2.7} has, for each $n$,  the solution$$Z_{n+1}=e^Wz_{n+1},$$where $z_{n+1}$ is the solution to the random differential equation
\begin{equation}\label{e2.8}
\barr{l}
\dd\frac d{dt}\,z_{n+1}+\frac1\lbb\,e^{-W}J(e^Wz_{n+1})+(\mu +\nu) z_{n+1}=e^{-W}F_n,\vsp
z_{n+1}(0)=x.\earr\end{equation}

If $F:L^2((0,T)\times\ooo;V')\to L^2((0,T)\times\ooo;V)$ is the linear continuous operator defined by
$$F(f)=Y,$$
where $Y$ is the solution to the stochastic equation
$$dY+\frac1\lbb\,J(Y)dt+\nu Y\,dt= Y\,dW+f\,dt;\ Y(0)=x,$$then we may rewrite \eqref{e2.2}--\eqref{e2.3} as
$$X_{n+1}=(\lbb(A-\nu I)+J)^{-1}[JF((\lbb(\nu I-A)+J)(X_n))+\lbb(A-\nu I)X_n],\ n=0,1,...$$
Equivalently,
\begin{equation}\label{e2.9}
X_{n+1}=\Gamma^n X_0,\ \ff n\in\nn,\end{equation}
where $\G: L^2((0,T)\times\ooo;V)\to L^2((0,T)\times\ooo;V)$ is the Lipschitzian and given by
\begin{equation}\label{e2.10}
\G=(\lbb(A-\nu I)+J)^{-1}[JF(\lbb(\nu I-A)+J)+\lbb(A-\nu I)].\end{equation}
Then, by Theorem \ref{t2.1}, we get

\begin{corollary}\label{c2.1} Under assumptions {\rm(i)-(iv)}, \eqref{e2.1}, for each $\lbb>0$ the solution $X$ to \eqref{e1.1} is expressed as
\begin{equation}\label{e2.11}
X=w-\lim_{n\to\9}\G^n X_0\mbox{\ \ in }L^2((0,T)\times\ooo;V),\end{equation}
where $X_0\in L^2((0,T)\times\ooo;V)$ is an arbitrary $(\calf_t)_{t\ge0}$-adapted process.
\end{corollary}
Here $w-\lim$ indicates the weak limit.

\section{Proof of Theorem \ref{t2.1}}
\setcounter{theorem}{0}
\setcounter{equation}{0}

Proceeding as in \cite{2}, we consider the spaces $\calh,\calv$ and $\calv'$, defined as follows. $\calh$ is the Hilbert space of all $(\calf_t)_{t\ge0}$-adapted processes $y:[0,T]\to H$ such that
$$|y|_\calh=\(\E\int^Y_0|e^{W(t)}y(t)|^2_Hdt\)^{\frac12}<\9,$$where $\E$ denotes the expectation in the  probability space $(\ooo,\calf,\mathbb{P})$. The space $\calh$ is endowed with the norm $|\cdot|_\calh$ generated by the scalar product
$$\<y,z\>_H=\E\int^T_0\<e^{W(t)}y(t),e^{E(t)}y(t)\>dt.$$
$\calv$ is the space of all $(\calf_t)_{t\ge0}$-adapted processes $y:[0,T]\to V$ such that
$$|y|_\calv=\(\E\int^T_0|e^{W(t)}y(t)|^2_Vdt\)^{\frac12}<\9.$$
$\calv'$ (the dual of $\calv)$ is the space of all $(\calf_t)_{t\ge0}$-adapted processes $y:[0,T]\to V'$ such that
$$|y|_{\calv'}=\(\E\int^T_0|e^{W(t)}y(t)|^2_{V'}dt\)^{\frac12}<\9.$$We have $\calv\subset\calh\subset\calv'$ with continuous  and dense embeddings. Moreover, $$_{\calv'}\<u,v\>_\calv=\E\int^T_0\<e^{W(t)}u(t),e^{W(t)}v(t)\>dt,\ v\in\calv,\ u\in\calv',$$is the duality pairing between $\calv$ and $\calv'$, with the pivot space $\calh$, that is,
$$_{\calv'}\<u,v\>_\calv=\<u,v\>_\calh,\ \ff u\in\calh,\ v\in \calv.$$
Now, for $x\in H$, define the operators $\cala:\calv\to\calv'$ and $\calb:D(\calb)\subset\calv\to\calv'$ as follows:
\begin{equation}\label{e3.1}
\barr{lcl}
(\cala y)(t)&=&e^{-W(t)}A(t)(e^{W(t)}y(t))-\nu y(t),\ \mbox{ a.e. }t\in(0,T),\ y\in\calv,\vsp
(\calb y)(t)&=&\dd\frac{dy}{dt}\,(t)+(\mu+\nu)y(t),\mbox{ a.e. }t\in(0,T),\ y\in D(\calb),\earr\end{equation}
\begin{equation}\label{e3.2}
\barr{r}
D(\calb)=\Bigl\{y\in\calv:y\in AC([0,T];V')\cap C([0,T];H),\ \pas,\\ \dd\frac{dy}{dt}\in\calv',\ y(0)=x\Bigr\}.\earr\end{equation}
Here, $AC([0,T];V')$ is the space of all absolutely continuous $V'$-valued functions on $[0,T]$. If $y\in D(\calb)$, then $y\in C([0,T];H)$ and $\frac{dy}{dt}$ is the derivative of $y$ in the sense of $V'$-valued distributions on $(0,T)$.
Then, equation \eqref{e1.14} can be expressed as
\begin{equation}\label{e3.3}
\calb y+\cala y=0.\end{equation}
Then, the map  $\wedge:\calv\to\calv'$ defined by
\begin{equation}\label{e3.4}
\wedge v=e^{-W}J(e^Wv),\ v\in V,\end{equation}is the canonical isomorphism of $\calv$ onto $\calv'$ and the scalar product $_\calv\<\cdot,\cdot\>_\calv$ of the space $\calv$ can be expressed as
\begin{equation}\label{e3.5}
_\calv\<v,\bar v\>_\calv={}_\calv\<v,\wedge\bar v\>_{\calv'},\ \ff v,\bar v\in\calv.\end{equation}We set
\begin{eqnarray}
(\cala^*u)(t)&\!\!\!=\!\!\!&\wedge\1\cala u(t)=e^{-W}J\1(A(t)(e^Wu)-\nu e^Wu),\, \ff u\in\calv,\ \ \ \ \ \label{e3.6}\\[2mm]
(\calb^*u)(t)&\!\!\!=\!\!\!&\wedge\1\calb u(t)=e^{-W}J\1\(e^W\(\dd\frac{du}{dt}\,+(\mu+\nu)u\)\),\label{e3.7}\\
 &&\hspace*{5cm}\ff u\in D(\calb^*)=D(\calb).\nonumber\end{eqnarray}
Since the operators $\cala$, $\calb$ and $\cala+\calb$ are maximal monotone in $\calv\times\calv'$ (\cite{2}, Lemma 4.1, Lemma 4.2), it is easily seen by \eqref{e3.6}-\eqref{e3.7} that $\cala^*$, $\calb^*$ and $\cala^*+\calb^*$ are maximal monotone in $\calv\times\calv$.

On the  other hand, by \eqref{e3.3} we can rewrite equation \eqref{e3.3} as
\begin{equation}\label{e3.8}
\calb^*y+\cala^*y=0.\end{equation}
Let $y\in D(\calb)$ be the unique solution to equation \eqref{e3.3} (see \cite{2}, Proposition 3.3). Then, $y$ is also the solution to  \eqref{e3.8}  and so, by Theorem 1 in \cite{6} (see, also, Corollary 6.1 in \cite{5}), we have that
\begin{equation}\label{e3.9}
y=\lim_{n\to\9}(I+\lbb\cala^*)\1 v_n\mbox{ weakly in }\calv\mbox{ as }n\to\9,\end{equation}where $\{v_n\}\subset\calv$ is, for $n\ge0$,  defined by
\begin{equation}\label{e3.10}
v_{n+1}=(I+\lbb\calb^*)\1(2(I+\lbb\cala^*)\1v_n-v_n)+
(I-(I+\lbb\cala^*)\1)v_n,\end{equation}
and $v_0$ is arbitrary in $\calv$. Here, $I$ is the identity operator in $\calv$.

The splitting algorithm \eqref{e3.9}--\eqref{e3.10} is just the Douglas--Rachford algorithm (\cite{4}) for equation \eqref{e3.8} and it can be equivalently expressed as
\begin{align}
y&=\lim_{n\to\9}y_n\mbox{ weakly in }\calv,\label{e3.11}\\
y_n&=(I+\lbb\cala^*)\1 v_n,\ n=0,1,...,\label{e3.12}\\
y_{n+1}+\lbb\cala^*y_{n+1}&=z_{n+1}+v_n-y_n,\label{e3.13}\\
z_{n+1}+\lbb\calb^*z_{n+1}&=2y_n-v_n,\label{e3.14}
\end{align}
where $v_0\in\calv$. (To get \eqref{e3.12}-\eqref{e3.14} from \eqref{e3.10}, we have used the identity\break $(I+\lbb\calb^*)\1(v+\lbb\calb^*v)=v,\ \ff v\in D(\calb^*) $ and the linearity of $\calb^*$.)

In fact, the weak convergence of $\{v_n\}$ in the space $\calv$ is also a consequence of the convergence of the Rockafellar proximal point algorithm \cite{7} for the maximal monotone operator $v\to G\1(v)-v$, where
\begin{equation}\label{e3.15}
G(z)=(I+\lbb\calb^*)\1(2(I+\lbb\cala^*)\1z-z)+z-(I+\lbb\cala^*)\1z,\ \ff z\in\calv.\end{equation}
(See \cite{5}, Theorem 4.)
Taking into account \eqref{e3.6}, \eqref{e3.7}, \eqref{e3.12} we rewrite \eqref{e3.14} as
\begin{equation}\label{e3.16}
\barr{l}
e^{-W}J(e^Wz_{n+1})+\lbb\(\dd\frac{dz_{n+1}}{dt}+(\mu+\nu)z_{n+1}\)\vsp
\qquad=e^{-W}J(e^W(2y_{n}-v_n))
=e^{-W}J(e^W(-\lbb\cala^* y_n+y_n))\vsp
\qquad=-\lbb e^{-W}A(t)(e^Wy_n)+\lbb\nu y_n+e^{-W}J(e^Wy_n)\earr
\end{equation}
and \eqref{e3.13} as
\begin{equation}\label{e3.17}
\barr{r}
 J(e^Wy_{n+1})+\lbb  A(t)(e^Wy_{n+1})-\lbb\nu e^W y_{n+1}\vsp
= J(e^W(z_{n+1}+v_n-y_n)).\earr\end{equation}
We set
$$X_n=e^Wy_n,\ Z_n=e^Wz_n   .$$Then, by \eqref{e3.16}, we get via It\^o's formula (see \cite{2} and \eqref{e2.8}, \eqref{e2.7})
$$\barr{l}
\lbb dZ_{n+1}+J(Z_{n+1})dt+\lbb\nu Z_{n+1}dt
=\lbb Z_{n+1}dW-\lbb A(t)X_ndt\vsp\hfill+\lbb\nu X_ndt+J(X_n)dt,\\
Z_{n+1}(0)=x.\earr$$
By \eqref{e3.17} and \eqref{e3.12}, we also get  that
$$\barr{l}
\lbb A(t)X_{n+1}(t)+J(X_{n+1}(t))-\lbb\nu X_{n+1}(t)\vsp
\qquad\quad=J(Z_{n+1}(t))+\lbb A(t)X_n(t)-\lbb\nu X_n(t),\ t\in(0,T),\earr$$which are just equations \eqref{e2.2}, \eqref{e2.3}. Moreover, by \eqref{e3.11}, we see that \eqref{e2.5} holds.

Assume now that $A(t):V\to V'$ is odd. Then so is $\cala^*:\calv\to\calv$ and also the operator $G$ defined by \eqref{e3.15}. Then, according to a result of J. Baillon \cite{1a}, the sequence $\{v_n\}$ defined by \eqref{e3.10}, that is $v_{n+1}=G(v_n)$, is strongly convergent in $\calv$. Recalling \eqref{e3.9}, we infer that so is the sequence $\{y_n\}$ and, consequently, \eqref{e2.6} holds. This completes the proof of  Theorem \ref{t2.1}.

\begin{remark}\label{r2.1}\rm One might expect that a similar splitting scheme can be constructed for nonlinear monotone operators $A(t):V\to V'$, where $V$ is a reflexive Banach space and $A(t)$ are demicontinuous coercive and with polynomial growth as in \cite{2}.  In fact, in this case, one might replace \eqref{e2.2} by
$$\barr{l}
\lbb dZ_{n+1}+Z_{n+1}dt+\lbb\nu Z_{n+1}dt=Z_{n+1}dW
-\lbb A_H(t)X_ndt
+\lbb\nu X_ndt+X_ndt,\\\hfill t\in(0,T),\vsp
\lbb A_H(t)X_{n+1}+X_{n+1}-\lbb\nu X_{n+1}=Z_{n+1}+\lbb A_H(t)X_n-\lbb\nu X_n,\earr$$
where $A_H(t)u=A(t)u\cap H.$ This question will be addressed in Section 5 below (see Remark \ref{r5.2}).\end{remark}

\section{Examples}
\setcounter{theorem}{0}
\setcounter{equation}{0}

We shall illustrate here the  splitting algorithm \eqref{e2.2}--\eqref{e2.3} for a few parabolic stochastic differential equations.

\begin{example}{\it Nonlinear stochastic parabolic equations.}\label{ex4.1}\mk

   \n\rm Consider the reaction-diffusion stochastic equation in $\calo\subset\rr^d$,
\begin{equation}\label{e4.1}
\barr{l}
dX-{\rm div}(a(t,\xi,\nabla X))dt+\nu X dt+\psi(X)dt=X dW\mbox{ in }(0,T)\times\calo,\vsp
X=0\mbox{ on }(0,T)\times\pp\calo,\ \ \ X(0)=x\mbox{ in }\calo.\earr\end{equation}
Here, $a:(0,T)\times\calo\times\rr^d\to\rr^d$ is measurable in $(t,\xi,r)$ continuous in $r$ on $\rr^d$,  $a(t,\xi,0)=0$. (The more general case, when $a:(0,T)\times\calo\times\ooo\times\rr^d)\to\rr^d$ is progressively measurable, could also be considered.) We assume also that
$$\barr{ll}
(a(t,\xi,r_1)-a(t,\xi,r_2))\cdot(r_1-r_2)\ge0,&\ff r_1,r_2\in\rr^d,\ (t,\xi)\in(0,T)\times\calo,\vsp
a(t,\xi,r)\cdot r\ge a_1|r|^2_d+a_2,&\ff r\in\rr^d,\ (t,\xi)\in(0,T)\times\calo,\vsp
|a(t,\xi,r)|_d\le c_1|r|_d+c_2,&\ff r\in\rr^d,\ (t,\xi)\in(0,T)\times\calo,\earr$$
where $a_1,c_1,\nu>0,\ a_2,c_2\in\rr,$ are independent of $(t,\xi)$, and $\psi:\rr\to\rr$  is a continuous and mo\-no\-to\-ni\-cally nondecreasing function such that $\psi(0)=0$ and $|\psi(r)|\le C(|r|^{\frac{2d}{d+2}}+1),$ $\ff r\in\rr$. Here $\calo\subset\rr^d$ is a bounded open subset with smooth boundary $\pp\calo$, and $|\cdot|_d$ is the Euclidean norm of $\rr^d$.

If $H=L^2(\calo)$, $V=H^1_0(\calo),$ $V'=H\1(\calo)$ and , for $t\in(0,T)$, the operator $A(t):V\to V'$ is defined by
$$_{V'}\<A(t)y,\vf\>_V=\int_\calo(a(t,\xi,\nabla y)\cdot\nabla\vf+\psi(y)\vf)d\xi,\ \ff\vf\in H^1_0(\calo),\ y\in H^1_0(\calo),$$then   Hypotheses (i)--(iii) are satisfied. As regards the Wiener process $W$, we assume here that, besides \eqref{e1.7}, the following condition holds:
$$\sum^\9_{j=1}\mu^2_j|\nabla e_j|^2_\9<\9.$$Then, by Theorem \ref{t2.1}, where $H$, $V$ and $A(t)$ are defined above and $J=-\Delta$ with Dirichlet homogeneous boundary conditions, if $x\in H^1_0(\calo)$, the solution $$X\in L^2(\ooo;C([0,T];L^2(\calo))\cap L^2((0,T)\times\ooo;H^1_0(\calo)))$$  to \eqref{e4.1} can be obtained as
\begin{equation}\label{e4.2}
X=w-\lim_{n\to\9} X_n\mbox{\ \ in }L^2((0,T)\times\ooo;H^1_0(\calo)),\end{equation}where $(X_n,Z_n)\in L^2((0,T)\times\ooo;H^1_0(\calo))$ is the  solution to the system (we~take $\lbb=1$)
\begin{equation}\label{e4.3}
\barr{l}
dZ_{n+1}-\D Z_{n+1}dt+\nu Z_{n+1}dt= Z_{n+1}dW+{\rm div}(a(t,\xi,\nabla X_n))dt-\D X_ndt\vsp
\hfill\mbox{ in }(0,T)\times\calo,\vsp
Z_{n+1}(0)=x\ \ \mbox{in }\calo,\vsp
Z_{n+1}=0\ \ \mbox{ in }(0,T)\times\pp\calo,\vsp
{\rm div}\,a(\nabla X_{n+1})+\D X_{n+1}=\D Z_{n+1}+{\rm div}(a(t,\xi,\nabla X_n))\mbox{ in  }(0,T)\times\calo,\vsp
\earr\end{equation}
where $X_0\in L^2((0,T)\times\ooo;H^1_0(\calo))$ is arbitrary but $\calf_t$-adapted.
Moreover,  if $a(t,\xi,-r)\equiv-a(t,\xi,r),\ \ff r\in\rr^d$, then the convergence \eqref{e4.2} is strong in $L^2((0,T)\times\ooo;H^1_0(\calo)).$
\end{example}

\begin{example}\label{ex4.2} {\it Stochastic porous media equations.}\mk

\n Consider the stochastic equation
\begin{equation}\label{e4.4}
\barr{ll}
dX-\D\psi(t,\xi,X)dt-\nu\D X dt=XdW&\mbox{ in }(0,T)\times\calo,\vsp
X(0,\xi)=x(\xi)&\mbox{ in }\calo,\vsp
\psi(t,\xi,X(t,\xi))=0&\mbox{ on }(0,T)\times\pp\calo,\earr \end{equation}
where $\calo$ is a bounded domain in $\rr^d$, $\nu>0$, the function \mbox{$\psi:[0,T]\times\ov\calo\times\rr\to\rr$} is continuous, $r\to\psi(t,\xi,r)$ is monotonically increasing in $r$, and there exist $a\in(0,\9)$ and $c\in[0,\9)$ such that
\begin{equation}\label{e4.5}
\barr{lcll}
r\psi(t,\xi,r)&\ge& a|r|^2-c,&\ff r\in\rr,\ (t,\xi,r)\in[0,T]\times\ov\calo,\vsp
|\psi(t,\xi,r)|&\le&c(1+|r|),&\ff r\in\rr,\ (t,\xi,r)\in[0,T]\times\ov\calo.\earr\end{equation}
We shall write equation \eqref{e4.4} under the form \eqref{e1.1} with $H=H^{-1}(\calo)$. Namely, we take $V=L^2(\calo),\ H=H\1(\calo),$ and $V'$ is the dual of $V$ with the pivot space $H\1(\calo)$. Then, $V\subset H\subset V'$ and
$$V'=\{\theta\in \cald'(\calo):\theta=-\D v,  \ v\in L^2(\calo)\},$$where $\D$ is taken in the sense of distributions on $\calo$. (Here $\cald'(\calo)$ is the space of Schwartz distributions on $\calo.$) The duality $_{V'}\<\cdot,\cdot\>_V$ is defined~as
$$_{V'}\<\theta,u\>_V=\int_\calo\wt\theta u\,d\xi,\ \ \wt\theta=(-\D)\1\theta,$$where
$\D$ is the Laplace operator with homogeneous Dirichlet boundary conditions on $\pp\calo$. The duality mapping $J:V\to V'$ is just the operator $-\Delta$ defined from $L^2(\calo)$ to $V'\subset\cald'(\calo)$ by
 $$\Delta u(\vf)=\int_\calo u\Delta\vf\,d\xi,\ \ff\vf\in H^1_0(\calo)\cap H^2(\calo).$$

The operator $A(t):V\to V'$ is defined by
$$_{V'}\<A(t)y,v\>_V=\int_\calo\psi(t,\xi,y)v\,d\xi,\ \ \ff y,v\in V=L^2(\calo),\ t\in[0,T].$$
Then, Hypotheses (i)--(iv) hold and so, if $x\in L^2(\calo)$, by Theorem \ref{t2.1}, the solution $X\in L^2(\ooo;C([0,T];H^{-1}(\calo))\cap L^2((0,T)\times\ooo;L^2(\calo)))$ to \eqref{e4.4} is given~by
$$\dd X=w-\lim_{n\to\9} X_n\mbox{\ \ in }L^2((0,T)\times\ooo;L^2(\calo)),$$
where
\begin{eqnarray}
&dZ_{n+1}-\Delta Z_{n+1}dt+\nu Z_{n+1}dt=Z_{n+1}dW-\D\psi(t,\cdot,X_n)dt
-\Delta Z_ndt\nonumber\\
&\hfill\mbox{ in }(0,T)\times\calo,\nonumber\\[2mm]
&Z_{n+1}(0)=x\in L^2(\calo),\ n=0,1,...,\label{e4.6}\\[2mm]
&\D\psi(t,\cdot, X_{n+1})+X_{n+1}=Z_{n+1}+\D\psi(t,\cdot, X_n),\ \mbox{ in }\calo,\nonumber\\
&\psi(t,\cdot,X_{n+1}(t,\cdot))=0\ \ \mbox{on }\pp\calo,\nonumber\\
&n=0,1,...,
 X_0\in L^2(0,T;L^2(\ooo;L^2(\calo))).\nonumber
\end{eqnarray}

If $\psi(t,\xi,r)=-\psi(t,\xi,-r),\ \ff r\in\rr,$ then the convergence of the sequence $\{X_n\}$ is strong in $L^2((0,T)\times\ooo;L^2(\calo)).$
\end{example}

\section{The case where $A(t)$ is maximal monotone in $H\times H$}
\setcounter{equation}{0}

Consider now equation \eqref{e1.1} under the following assumptions on $A$:
\bit\item[(j)] {\it$A:[0,T]\times H\times\ooo\to H$ is progressively measurable and, for each $(t,\oo)\times[0,T]\times\ooo$ the operator $u\to A(t,\oo,u)$ is maximal monotone in $H\times H$. Moreover, there is $f\in L^2((0,T)\times\ooo;H)$ such that
\begin{equation}\label{e5.1}
(I+A(t))^{-1}f(t)\in L^2((0,T)\times\ooo;H).\end{equation}}\eit
We assume also that condition \eqref{e2.1} holds.

It should be noted that, if $A(t):V\to V'$ satisfies assumptions (i)-(ii), where $V$ is a reflexive Banach space, then the operator $A(t):H\to H$, defined~by
$$A(t)_Hu=A(t)u\cap V,$$ satisfies assumption (j). However, the class of the operators $A$ satisfying (j) is considerably larger.

We consider the splitting scheme (which is well defined by strong monotonicity of $\cala^*_1+\calb^*_1$)
\begin{equation}\label{e5.2}
\barr{l}
\lbb dY_{n+1}+(1+\lbb\nu)Y_{n+1}dt=\lbb Y_{n+1}dW\vspp\qquad\  \quad+\,(V_n-((1-\lbb\nu)I-\lbb A(t))^{-1}V_n)dt,\vsp
Y_{n+1}(0)=x\mbox{ in }(0,T),\vsp
V_{n+1}=Y_{n+1}+V_n-((1-\lbb\nu)I-\lbb A(t))^{-1}V_n,\earr
\end{equation}
where $V_0\in L^2((0,T)\times\ooo;H)$ is an $(\calf_t)_{t\ge0}$-adapted process such that $A(t)V_0\in L^2((0,T)\times\ooo;H)$. We have

\begin{theorem}\label{t5.1} Assume that $x\in H$ and that equation \eqref{e1.1} has a solution $X\in L^2(\ooo;C([0,T];H))$ such that $A(t)X\in L^2((0,T)\times\ooo;H)$.

Then, for $n\to\9$,
\begin{equation}\label{e5.3}
V_n\to V\mbox{ weakly in }L^2((0,T)\times\ooo;H),
\end{equation}
where $X=((1-\lbb\nu)I+A(t))^{-1}V$ is the solution to \eqref{e1.1}.

If $A(t)$ is odd, then the convergence \eqref{e5.3} is strong.
\end{theorem}

\pf The operators $\cala^*_1$ and $\calb_1^*$ defined by
$$\barr{rcll}
(\cala^*_1u)(t)&=&e^{-W}A(t)(e^Wu)-\nu u,&\ff u\in D(\cala^*_1),\vsp
(\calb^*_1u)(t)&=&\dd\frac{du}{dt}+(\mu+\nu)u,&\ff u\in D(\calb^*_1),\earr$$
with the domains
$$\barr{rcll}
D(\cala^*_1)&=&\{u\in\calh;\ e^{-W}A(t)(e^Wu)-\nu u\in\calh\},\vsp
D(\calb^*_1)&=&\{u\in\calh;\ u\in W^{1,2}([0,T];H)\cdot\pas,\ u(0)=x\}\earr$$
are, by the above hypotheses, maximal monotone in $\calh\times\calh$ (see also \cite{3}). Moreover, there is at least one solution $y^*$ to the equation
\begin{equation}\label{e5.4}
\cala^*_1 y^*+\calb_1^*y^*=0.
\end{equation}
Then, again by \cite{6}, it follows that the sequence $\{v_n\}\subset\calh$ defined by
\begin{equation}\label{e5.5}
\barr{r}
v_{n+1}=(I+\lbb\calb^*_1)^{-1}
(2(I+\lbb\cala^*_1)^{-1}v_n-v_n)
+v_n-(I+\lbb\cala^*_1)^{-1}v_n,\\ n=0,1,...\earr
\end{equation}
is weakly convergent in $\calh$ to $v^*$, where $(1+\lbb\cala^*_1)^{-1}v^*=y^*$ is the solution to equation \eqref{e5.4}.

We set
\begin{equation}\label{e5.6}
\wt z_{n+1}=v_{n+1}-v_n+(I+\lbb\cala^*_1)^{-1}v_n
\end{equation}
and, by \eqref{e5.5}, we have
\begin{equation}\label{e5.7}
\wt z_{n+1}+\lbb \calb^*_1 z_{n+1}=v_n-(I+\lbb\cala^*_1)^{-1}v_n.
\end{equation}
Then, if   $Y_n=e^W\wt z_n$ and $V_n=e^Wv_n$, we can rewrite \eqref{e5.6}-\eqref{e5.7} as  \eqref{e5.2} and get \eqref{e5.3}, as claimed.

\begin{remark}\label{r5.2}\rm The convergence of the splitting algorithm \eqref{e5.1}-\eqref{e5.2} does not require conditions of the form (ii)-(iii) for the operator $A(t)$ but in change it requires the existence of a sufficiently regular solution $X$ for equation \eqref{e1.1} $(A(t)x\in L^2((0,T)\times\ooo;H)$) which is not the case for Examples \ref{ex4.1}, \ref{ex4.2}. Such a condition holds, however, for the stochastic reaction-diffusion equation
$$\barr{l}
dX-\Delta X\,dt+\Psi(X)dt=X\,dW\mbox{ in }(0,T)\times\calo,\vsp
X=0\mbox{ on }(0,T)\times\pp\calo,\vsp
X(0)=x,\earr$$
if $x\in H^1_0(\calo)$ and $\Psi:\rr\to\rr$ is continuous and monotonically increasing and for other stochastic parabolic equations as well. \end{remark}

\bk\n{\bf Acknowledgments.} Financial support through SFB~701 at Bielefeld University is gratefully acknowledged.


\begin{thebibliography}{nn}

\bibitem{1a} J.B. Baillon, Quelques propri\'et\'es de convergence asymptotiques pour les contractions impaires, {\it C.R. Acad. Sci. Paris}, A, 283 (1976), 587-590.
\bibitem{1}V. Barbu, {\it Nonlinear Differential Equations of Monotone Type in Banach Spaces}, Springer, New York, 2010.

\bibitem{2} V. Barbu, M. R\"ockner, An operatorial approach to stochastic partial differential equations driven by linear multiplicative noise, {\it J. Eur. Math. Soc.}, 17 (2015), 1789-1815.

\bibitem{3} H. Brezis, P.L. Lions, Produits infinis de resolvantes, {\it Israel J. Math.}, 29 (1978), 329-345.


\bibitem{4a} G. Da Prato, J. Zabczyk, {\it Stochastic Equations in Infinite Dimensions}, Encyclopedia of Mathematics and Its Applications, Cambridge University Press, Cambridge, 2014.

\bibitem{4} J. Douglas, H.H. Rachford, On the numerical solution of heat conduction problems in two and three space variables, {\it Trans. Amer. Math. Soc.}, 82 (1956), 421-439.



\bibitem{5} J. Eckstein, D.P. Bertsekas, On the Douglas-Rachford splitting method and the proximal point algorithm for maximal monotone operators, {\it Math. Programming}, 55 (1992), 293-318.



\bibitem{6} P.L. Lions, B. Mercier, Splitting algorithm for the sum of two nonlinear operators, {\it SIAM J. Numerical Analysis}, 16 (1979), 277-293.

\bibitem{7} R.T. Rockafellar, Monotone operators and the proximal point algorithm, {\it SIAM J. Control Optimization}, 14 (1976), 877-898.

\end{thebibliography}
\end{document}